\documentclass{amsart}
\usepackage{amssymb}
\usepackage[all]{xy}
\CompileMatrices
\newtheorem*{introthm}{Theorem}
\newtheorem{thm}{Theorem}[section]
\newtheorem{lemma}[thm]{Lemma}

\newtheorem{coro}[thm]{Corollary}
\theoremstyle{definition}
\newtheorem{defi}[thm]{Definition}
\newtheorem{exam}[thm]{Example}
\newtheorem{rem}[thm]{Remark}

\def\div{{\rm div}}

\def\reg{{\rm reg}}
\def\tq{\mathbin{\!{/ \! \! \lower 3pt \hbox{$\scriptscriptstyle{\rm
          tq}$}}\!}}
\def\rq#1{\widehat{#1}}
\def\t#1{\widetilde{#1}}

\def\KK{{\mathbb K}}
\def\TT{{\mathbb T}}

\def\QQ{{\mathbb Q}}

\def\Of{\mathcal{O}}
\def\id{{\rm id}}
\def\Div{\operatorname{Div}}
\def\WD{\operatorname{Div}_{W}}
\def\CaD{\operatorname{Div}_{C}}
\def\WDi{\WD^{\TT}}

\def\PDiv{\operatorname{PDiv}}
\def\PD{\operatorname{PDiv}}

\def\ClDiv{\operatorname{ClDiv}}

\renewcommand{\div}{\operatorname{div}}
\renewcommand{\phi}{\varphi}
\def\Hom{{\rm Hom}}
\def\codim{{\rm codim}}

\def\Spec{{\rm Spec}}

\def\pr{{\rm pr}}

\def\lra{\longrightarrow}
\newcounter{itemnumber}
\newenvironment{itemlist}{%
\begin{list}{(\roman{itemnumber})}{\usecounter{itemnumber}%
\topsep0pt\itemsep1ex}}{\nobreak\end{list}\setcounter{itemnumber}{1}}
\begin{document}
\title{Lifting of morphisms to quotient presentations}
\author{Florian Berchtold}
\address{Fachbereich Mathematik und Statistik, Universit\"at Konstanz,
78457 Konstanz, Germany}
\email{florian.berchtold@uni-konstanz.de}
\begin{abstract}\noindent
In this article we investigate algebraic morphisms between toric
varieties. Given presentations of toric varieties as quotients we are
interested in the question when a morphism admits a lifting to these quotient
presentations. We show that this can be completely answered in terms of
invariant divisors. As an application we prove that two toric varieties, which
are isomorphic as abstract algebraic varieties, are even isomorphic as toric
varieties. This generalizes a well-known result of Demushkin on affine toric
varieties.
\end{abstract}
\maketitle
\section*{Introduction}
By definition, a toric variety is a normal algebraic variety $X$ together with
an effective action of an algebraic torus $\mathbb{T}$ which imbeds
equivariantly as a dense open subset in $X$. Generalising projective spaces
toric varieties admit presentations as quotients of quasi-affine toric
varieties. There are different possibilities for such presentations, see
\cite{Cox1}, \cite{BV} and \cite{Ka}, for example. A systematic treatment of
all possible quotient presentations is given in \cite{Ajs}.

It is a classical fact that in the case of projective space each morphism
admits a lifting to the affine spaces above. We generalize this result to
smooth toric varieties. For singular toric varieties we give a necessary and
sufficient criterion when a morphism admits a lifting.

The main result (Theorem \ref{mainthm}) is a classification of all possible
liftings of a given morphism for fixed quotient presentations in terms of
homomorphisms defined on the groups of invariant divisors. As an application
we prove

\begin{introthm}
Let $X$ and $X'$ be toric varieties which are isomorphic as abstract
algebraic varieties. Then they are also isomorphic in the category of toric
varieties.
\end{introthm}

For smooth complete toric varieties this is an immediate consequence from the
fact that in this case the automorphism group is linear algebraic with maximal
torus $\TT$ (compare Demazure \cite{Dem}) and in such a group any two maximal
tori are conjugate according to Borel's theorem (\cite{Bo}). For affine toric
varieties the proof is much harder and is due to Demushkin \cite{De} and
Gubeladze \cite{Gu}. Our approach is completely different and works for
arbitrary toric varieties. The main idea is to replace functions on the
varieties in question by the more general notion of Weil divisors (a similar
idea can be found in \cite{Fi}).

This article is divided into four sections. In the first section we mention
some basic facts concerning quotient presentations of toric varieties. In the
second section we present the main theorem and some immediate
consequences. The proof of this theorem will be given in the third
section. Finally in the fourth section we prove the toric isomorphism theorem.
%
%
\section{Quotient Presentations}
Let $\KK$ be an algebraically closed field. We consider non-degenerate toric
varieties over $\KK$. A toric variety $X$ is called non-degenerate, if $X$
does not admit an equivariant decomposition $X \cong X_{1}\times {\KK^{*}}$,
where $X_{1}$ is a toric variety of dimension $\dim X - 1$. Observe that a
variety $X$ is non degenerate if and only if for the ring of invertible
regular functions $\Of(X)^{*} = \KK^{*}$ holds. We often identify the torus
$\TT \cong {({\KK}^{*})}^{\dim X}$ with the unique open dense orbit in $X$.

Any toric variety admits a presentation as a quotient of a quasiaffine variety
by the action of a diagonalisable group. We briefly recall some basic facts
concerning this approach (cf. \cite{Cox1}, \cite{Ka}, \cite{Ajs}). Let a
surjective toric morphism $q \colon \rq{X} \to X$ be given. The strict
transform of the group of Weil divisors $\WD(X)$ is defined by the composition
\[ q^{\sharp}\colon \WD(X) = \CaD(X_{\reg}) \overset{q^{*}}{\lra}
\CaD(q^{-1}X_{\reg}) \subset \WD(\rq{X}), \] where $\CaD(X)$ denotes the
subgroup of Cartier divisors and $X_{\reg}$ is the set of regular points of
$X$.
\begin{defi}\label{reprquot}
A surjective toric morphism $q \colon \rq{X} \to X$ is called a {\em quotient
presentation} of $X$, if the following conditions are fulfilled:
\begin{itemlist}
\item $\rq{X}$ is a quasi-affine toric variety,
\item the restriction of the strict transform $q^{{\sharp}}$ to the group
$\WDi(X)$ of invariant Weil divisors is bijective.
\end{itemlist}
\end{defi}
We now give a general construction to obtain all quotient presentations for a
given non-degenerate toric variety $X$. Let $M := \Hom(\TT, \KK^{*})$ denote
the character group of the torus $\TT$. Since $X$ is non-degenerate, $M$
injects into the group $\Div_{W}^{\TT}(X)$ of torus-invariant Weil divisors of
$X$ (cf. \cite{Ful}, \cite{Oda}). More precisely, it is isomorphic to the
subgroup of invariant principal divisors $\PDiv^{\TT}(X)$.

Let $\widehat{M}$ be a subgroup of $\Div_{W}^{\TT} (X)$ containing $M$,
i.e. $M \subset \widehat{M} \subset \Div_{W}^{\TT} (X)$. Assume further that
$\widehat{M}$ has enough divisors: that is, for each maximal $\TT$-stable
affine chart $U \subset X$ there is an effective divisor $D \in \widehat{M}$
such that $X \setminus U = |D|$ (cf. \cite{Ka}, \cite{Ajs}). Here $|D|$
denotes the support of the divisor $D$.

The idea to construct a quotient presentation is to make all divisors in
$\rq{M}$ principal. For that purpose let $S(\widehat{M})$ denote the algebra
$\KK[\widehat{M}_{\ge 0}]$ associated to the semigroup $\widehat{M}_{\ge 0}$
of effective invariant Weil divisors of $\widehat{M}$ and denote by $T^{D}$
the element of $S(\widehat{M})$ corresponding to the divisor $D \in
\rq{M}_{\ge 0}$. The finitely generated algebra $S(\widehat{M})$ defines an
affine toric variety $\Spec(S(\widehat{M})) =: \rq{U}$. Then we can view
$T^{D}$ as a regular function on $\rq{U}$. By $V(\rq{U}; T^{D_{1}}, \dots,
T^{D_{r}})$ we denote the common zero set of the functions $T^{D_{i}}$ on
$\rq{U}$. Set
\[ \rq{X} := \rq{U} \setminus Z(\widehat{M}), \text{
where } Z(\widehat{M}) := \bigcup_{|D_{1}| \cap \dots \cap |D_{r}| =
\emptyset} V(\rq{U}; T^{D_{1}}, \dots, T^{D_{r}})). \] Obviously, the
so-called exceptional set $Z (\widehat{M})$ has codimension at least two in
$\rq{U}$. Define the action of $H := \Hom ( \rq{M}/M, \KK^{*})$ by
\[ H \times S \to S, \; (h,T^{D}) \mapsto h\cdot T^{D} :=
(h(\overline{D}))^{-1} T^{D}, \] where $\overline{D}$ denotes the class of $D$
in $\rq{M}/M$. This action admits a good quotient $q \colon \rq{X} \to X$
(\cite{Ajs}, Proposition 3.2), that is an $H$-invariant, affine mapping such
that $q_{*}(\Of_{\rq{X}})^{H} = \Of_{X}$ holds.

Furthermore, any quotient presentation arises in that way:
\begin{thm}(\cite{Ajs}, Theorem 2.2)\label{equicat}
Any subgroup $\rq{M}$ of $\Div_{W}^{\TT}(X)$ containing $M$ with enough
divisors defines a quotient presentation of $X$. Conversely, given any
non-degenerate quotient presentations of $X$, there is a unique subgroup
$\rq{M}$ of $\Div_{W}^{\TT}(X)$ containing $M$ with enough divisors defining
this presentation.
\end{thm}

\begin{exam} Let $X$ be a toric variety.
\begin{itemlist}
\item (Cox's construction) The construction $\t{X}$ given in \cite{Cox1}
corresponds to $\rq{M} = \Div_{W}^{\TT}(X)$, the group of all invariant Weil
divisors. In this situation, the affine hull of $\t{X}$ is the affine space
$\KK^{\dim \t{X}}$.
\item (Kajiwara's construction) Assume that the group $\Div_{C}^{\TT} (X)$ of
all invariant Cartier divisors has enough divisors. Then, for $\rq{M} =
\Div_{C}^{\TT} (X)$ one obtains the quotient presentation $\rq{X}^{K}$ due to
Kajiwara (\cite{Ka}).
\end{itemlist}
\end{exam}
From Theorem \ref{equicat} it follows, that Cox's construction $\t{X}$ of the
last example is universal:
\begin{coro}(\cite{Ajs}, Example 2.3)\label{universal}
Let $q \colon \rq{X} \to X$ be any quotient presentation of $X$. Then the
quotient presentation $\pi \colon \t{X} \to X$ of Cox admits a unique
factorisation $\pi = q \circ \pi_{1}$, where $\pi_{1} \colon \t{X} \to \rq{X}$
is a quotient presentation of $\rq{X}$. In particular, applying Cox's
construction to $\rq{X}$ one obtains $\t{X}$.
\end{coro}
\section{Lifting of Morphisms}
Let $f\colon X' \to X$ be a morphism between non-degenerate toric varieties
and let $q' \colon \rq{X'} \to X'$ and $q \colon \rq{X} \to X$ be fixed
quotient presentations. Let $\widehat{M}'$ and $\widehat{M}$ denote the groups
of Weil divisors introduced in the last section corresponding to the quotient
presentations $q'$ and $q$, respectively. In this section we answer the
question, when there exists a morphism $\rq{F} \colon \rq{X'} \to \rq{X}$
making the following diagram commutative:
\begin{equation}\label{LP}
\xymatrix{%
{\rq{X'}} \ar[r]^{\rq{F}} \ar[d]_{q'} & {\rq{X}} \ar[d]^{q} \\ X' \ar[r]_{f} &
X} 
\end{equation}

In the case of projective space it is well known that the answer to this
question is always positive. Moreover each lifting is given by homogeneous
polynomials of the same degree. In our framework this translates to the
concept of equivariant morphisms: The groups $H' := \Hom (\widehat{M}'/M',
\KK^{*})$ and $H := \Hom (\widehat{M}/M, \KK^{*})$ act on $\rq{X'}$ and on
$\rq{X}$, respectively. A morphism $F \colon \rq{X'} \to \rq{X}$ is called
equivariant with respect to the actions of $H'$ and $H$, if there is a
homomorphism of algebraic groups $\psi \colon H' \to H$ such that $F (h'z) =
\psi(h') F(z)$ holds.

For any lifting $\rq{F}_{1}$ of $f$ and any $h \in H$ the morphism $\rq{F}_{2}
:= h \rq{F}_{1}$ is also a lifting. In the case of projective space any two
liftings of a given morphism differ only by such a constant. We shall call two
liftings $\rq{F}_{1}$ and $\rq{F}_{2}$ equivalent, if there is an $h \in H$
such that $\rq{F}_{1} = h \rq{F}_{2}$ holds.

It turns out that the existence of liftings is related to the existence of
certain pullbacks of Weil divisors:
\begin{defi}\label{compatible}
Let $f \colon X' \to X$ be a morphism. A homomorphism of groups $\phi \colon
\widehat{M} \to \Div_{W}(X')$ which extends the pullback of Cartier divisors
is called a {\em geometric pullback of $\widehat{M}$ with respect to $f$}, if
it fulfils the following two properties:
\begin{itemlist}
\item $\phi ({\widehat{M}}_{\ge 0}) \subset \Div_{W}(X')_{\ge 0}$,
\item $f(|\phi (D)|) \subset |D|$.
\end{itemlist}
\end{defi}
\begin{rem}\label{codim}
Let $f \colon X' \to X$ be a morphism with the property
\[ \codim_{X'} f^{-1} (X_{\text{sing}}) \ge 2. \]
Then the strict transform $f^{\sharp}$ of Weil divisors is well defined and
defines the unique geometric pullback with respect to $f$.
\end{rem}

\begin{thm}\label{mainthm} 
Let $X$ and $X'$ be non-degenerate toric varieties, and let $f\colon X' \to X$
be a morphism with $f(X') \cap \TT \neq \emptyset$. Let $q \colon \rq{X} \to
X$ and $q' \colon \rq{X'} \to X'$ be quotient presentations, corresponding to
the subgroups $\rq{M}$ and $\rq{M'}$. Then the equivalence classes of
equivariant liftings $\widehat{F} \colon \rq{X'} \to \rq{X}$ stay in a
one-to-one correspondence with the geometric pullbacks of $\widehat{M}$, the
images of which are contained in $\widehat{M'} + \PDiv(X')$.
\end{thm}
Before we give the proof of this result in the following section we firstly
present some corollaries, generalising results obtained by Cox \cite{Cox2} and
Kajiwara \cite{Ka}.
\begin{coro}\label{coxcon}
Let $f \colon X' \to X$ be a morphism with $f(X') \cap \TT \neq \emptyset$.
\begin{itemlist}
\item $f$ can be lifted to an equivariant morphism $\t{X'} \to \t{X}$, if and
only if there exists a geometric pullback of $\Div_{W}^{\TT} (X)$ with respect
to $f$.
\item If $\Div_{C}^{\TT}$ has enough divisors, then $f$ admits a unique
equivalence class of equivariant liftings $\t{X'} \to \rq{X}^{K}$. If
additionally $\Div_{C}^{\TT'}(X')$ has enough divisors, then there is even a
lifting $\rq{X'}^{K} \to \rq{X}^{K}$.
\end{itemlist}
\end{coro}

Using Remark \ref{codim} we immediately obtain
\begin{coro}\label{isolift}  
Let $\rq{X} \to X$ be a quotient presentation of $X'$, and let $f \colon X'
\to X$ be a morphism with $f(X') \cap \TT \neq \emptyset$. If $\codim_{X'}
f^{-1} (X_{\operatorname{sing}}) \ge 2$, then $f$ has a unique equivalence
class of equivariant liftings $\t{X'} \to \rq{X}$. In particular, if $X$ is
smooth, then there is precisely one equivalence class of equivariant liftings
$\t{X'} \to \rq{X}$ of any $f$.
\end{coro}
In the theory of toric varieties simplicial or $\QQ$-factorial toric varieties
play a special role. In this case it is easy to check whether a given pullback
is geometric:
\begin{rem}
Let $X$ be a simplicial toric variety, $\rq{M}$ a subgroup of $\Div_{W}^{\TT}
(X)$, and $f \colon X' \to X$ a morphism. A homomorphism of groups $\phi
\colon \widehat{M} \to \Div_{W}(X')$ is a geometric pullback, if and only if
it extends the usual pullback of Cartier divisors. In other words, the
conditions (i) and (ii) in \ref{compatible} are automatically fulfilled.
\end{rem}

In particular our results show that liftings do not always exist:
\begin{exam}
Let $X$ be the zero set of the function $z_{3}^{2} - z_{1}z_{2}$ in $\KK^{3}$
together with the torus embedding
$$ {(\KK^{*})^{2}} \to X; \; (t_{1},t_{2}) \mapsto
(t_{1}^{2}/t_{2},t_{2},t_{1}). $$ Consider the blow up $f \colon X' \to X$ of
the origin in $X$. There is no lifting $\t{X'} \to \t{X}$ of $f$ to the
constructions of Cox.
\end{exam}

\proof %
The $z_{1}$-axis which is contained in $X$ defines a torus-invariant prime
divisor $D_{1}$ in $X$. This is no principal divisor, but $2D_{1}$ is. Now,
$f^{*} (2D_{1}) = 2D'_{1} + E$, where $E$ denotes the exceptional divisor and
$D_{1}'$ is the preimage of $D_{1}$ outside of the origin. From this equation,
it is easy to see that the pullback $f^{*}$ of Cartier divisors does not admit
an extension to the group of invariant Weil divisors. By Corollary
\ref{coxcon}, a lifting cannot exist.  \endproof
\section{The Proof of Theorem \ref{mainthm}}
In this section we present the proof of Theorem \ref{mainthm}. Firstly we show
the easier part, namly how a lifting defines a geometric pullback. For that
purpose assume that $\rq{F} \colon \rq{X'} \to \rq{X}$ is an equivariant
lifting of $f \colon X' \to X$. Denote the acting groups on $\rq{X}$ and
$\rq{X'}$ by $H$ and $H'$, respectively.

For a divisor $D \in \rq{M}$ the rational function $T^{D}$ is mapped by
$\rq{F}^{*}$ to an $H'$-homogeneous rational function $r$, say of degree
$\overline{D'}$ for some $D' \in \rq{M'}$, because $\rq{F}$ is
equivariant. Then the rational function $r/T^{D'}$ is $H'$-invariant. Since $q'
\colon \rq{X'} \to X'$ is a good quotient, we conclude $r/T^{D'} = {q'}^{*}
(s)$ for some rational function $s$ on $X'$. Hence, we obtain that $r$
is of the form 
\begin{equation}\label{eq:decom}
r = T^{D'} {q'}^{*} (s), \quad D' \in \rq{M'}, \, s \in \mathcal{R}(X')
\end{equation}
Hence, defining $\phi \colon \rq{M} \to \rq{M'} + \PD(X')$ by $D \mapsto D' +
\div(s)$, we obtain a homomorphism of groups. Moreover we have
\begin{lemma}\label{pullbackprops}
The map $\phi$ is a geometric pullback.
\end{lemma}
\proof
Observe that for the pullback of Cartier divisors we have $\rq{F}^{*} q^{*} = 
{q'}^{*} f^{*}$, because $\rq{F}$ is a lifting. But by the definition of
$\phi$ we also have $\rq{F}^{*} q^{*} = {q'}^{*} \phi$. As ${q'}^{*}$ is
injective, we conclude that $\phi$ extends the pullback of Cartier divisors 

To show the first property of Definition \ref{compatible}, let $D \in \rq{M}$
be effective. Then the functions $T^{D}$ and $\rq{F}^{*} (T^{D})$ are
regular. But then $\phi(D)$ has to be an effective divisor, too.

The second property follows from 
\begin{align*}
|\phi (D)| & = q'(|{q'}^{\sharp} \phi (D)|) = q'(|\rq{F}^{*} q^{\sharp} (D)|)
\\ & \subset q' \rq{F}^{-1} q^{-1} (|D|) = q' {q'}^{-1} f^{-1} (|D|) = f^{-1}
(|D|),
\end{align*}
where the first equality holds, because good quotients map closed saturated
sets onto closed sets.  
\endproof

Thus we are left to show that a geometric pullback induces a lifting of $f$.
The main step of this part of the proof will be the construction of a lifting
$F \colon \t{X'}\to \rq{X}$, where $\pi' \colon \t{X'} \to X'$ again denotes
Cox's construction. From this we can easily deduce the existence of a lifting
$\rq{F} \colon \rq{X'} \to \rq{X}$. 

Given a morphism $f \colon X' \to X$ and a geometric pullback $\phi \colon
\rq{M} \to \rq{M'} + \PDiv(X')$ with respect to $f$, we shall construct the
lifting in such a way that for each $D \in \rq{M}$ we have
$$ \div F^{*} (T^{D}) = {\pi'}^{\sharp} \phi(D). $$ From this it will be clear
that both constructions are inverse to each other.

Let $B_{1}, \dots, B_{r}$ be a basis of $\rq{M}$. Consider the divisors
${\pi'}^{\sharp} \phi (B_{j})$ in $\t{X'}$ for $1 \le j \le r$ and recall that
the codimension of the exceptional set $Z(\Div_{W}^{\TT}(X'))$ in the affine
hull $\KK^{k'} = \Spec(S(\WDi(X')))$ of $\t{X'}$ is at least two. Therefore
for each $j$ the divisor ${\pi'}^{\sharp}\phi (B_{j})$ extends uniquely to a
divisor of $\KK^{k'}$. Hence it is given by a rational function $r_{j}$, which
is unique up to a constant. We now define for each $t \in {(\KK^{*})}^{r}$ a
homomorphism of the algebras of rational functions
\begin{equation}\label{eq:ratpb}
F_{t}^{*} \colon {\mathcal{R}}(\rq{X}) \to {\mathcal{R}} (\t{X'}) \text{ by }
T^{D_{j}} \mapsto t_{j}r_{j}.
\end{equation}
\begin{lemma}\label{comorprops}
The homomorphism $F^{*}_{t}$ has the following properties:
\begin{itemlist}
\item $F_{t}^{*} (S(\widehat{M})) \subset \mathcal{O} (\KK^{k'}) = \KK
[T'_{1},\dots, T'_{k'}]$,
\item For any $D \in \rq{M}$ there is an invariant divisor $D' \in
\widehat{M}'$ and a rational function $s$ on $X'$ such that $F_{t}^{*} (T^{D})
= {T}^{D'} (s\circ \pi')$ holds.
\end{itemlist}
\end{lemma}
\proof (i) If $T^{D}$ is in $S(\widehat{M})$, then $D$ is by definiton an
effective divisor. Thus from property (i) of Definition \ref{compatible} we
see that $\phi (D)$ is effective. But then ${\pi'}^{\sharp} \phi(D)$ is
effective, too. Therefore the rational function $r_{D} = F_{t}^{*} (T^{D})$ is
even a regular one as proposed.

(ii) Since $\phi (\widehat{M}) \subset \widehat{M}' + \PDiv(X')$, there is a
representation $\phi (D) = D' + \div(s)$ for an invariant divisor $D' \in
\widehat{M}'$ and a rational function $s$. By definition of $F^{*}_{t}$ and
linearity we have:
\[ \div (F_{t}^{*} (T^{D})) = {\pi'}^{\sharp} \phi (D) = {\pi'}^{\sharp}
(D' + \div(s)) = \div ({T}^{D'}) + \div (s\circ \pi'). \] Thus after replacing
$s$ by a suitable multiple we obtain the desired result.  
\endproof

From this Lemma we see that $F^{*}_{t}$ defines a morphism $\KK^{k'} \to
\rq{U} = \Spec(S(\rq{M}))$. By restriction we obtain a morphism $F_{t} \colon
\t{X'} \to \rq{U}$.
\begin{lemma}\label{morprops}
The morphism $F_{t}$ has the following properties:
\begin{itemlist}
\item $F_{t}(\t{X'}) \subset \rq{X}$.
\item The pullback $\phi$ induces a homomorphism
\[ \psi \colon G := \Hom(\ClDiv(X'),\KK^{*}) \to \Hom(\widehat{M}/M,
\KK^{*})=: H, \] such that $F_{t}$ is equivariant with respect to $\psi$.
\end{itemlist}
\end{lemma}

\proof (i) Assume the contrary, that is there is a $y \in \t{X'}$ which is
mapped into $Z(\widehat{M})$. Thus there are divisors $D_{1},\dots, D_{r} \in
\widehat{M}$ with an empty intersection, such that all functions $T^{D_{i}}$
vanish at $F_{t}(y)$. This implies $r_{D_{i}}(y) = 0$, i.e. $y \in |\div
(r_{D_{i}})|$, for all $i$. But by construction, $\div (r_{D_{i}}) =
{\pi'}^{\sharp} \phi (D_{i})$. Hence, from property (ii) in Definition
\ref{compatible} we obtain $f(\pi'(y)) \in |D_{i}|$ for all $i$. But this
contradicts the fact that the divisors $D_{1},\dots, D_{r}$ have an empty
intersection.

(ii) As $\phi \colon \widehat{M} \to \Div_{W}(X')$ is an extension of the
pullback of Cartier divisors, it obviously maps $M$ into
$\PDiv(X')$. Therefore it induces a homomorphism $[\phi] \colon \widehat{M}/M
\to \ClDiv(X')$, which in turn induces the dual homomorphism $\psi$. By Lemma
\ref{comorprops}(ii) the comorphism $F^{*}_{t}$ maps $T^{D}$ to the rational
function ${T}^{D'} (s\circ\pi')$, where $[D'] = [\phi] (\overline{D} )$. For a
$g \in G$ this yields
\begin{align*}
g \cdot F_{t}^{*} (T^{D}) & = g \cdot {T}^{D'} (s \circ \pi') = \bigl( g([\phi]
(\overline{D}) \bigr)^{-1} F_{t}^{*}(T^{D}) \\ & = \bigl( \psi (g)
(\overline{D}) \bigr)^{-1} F_{t}^{*}(T^{D}) = F_{t}^{*} (\psi(g) \cdot T^{D}).
\end{align*} 
Translating this into terms of the morphism $F_{t}$ gives that $F_{t}$ is in
fact equivariant with respect to $\psi$.
\endproof 

We now finally show that for a suitable choice of $t$ the map $F_{t}$ is
indeed a lifting of $f$. Firstly, by Lemma \ref{morprops}(i), $F_{t}$ defines
a morphism $\t{X'}\to \rq{X}$. Secondly, as $F_{t}$ is equivariant, it induces
a morphism $X' \to X$. Therefore it remains to show that this induced morphism
coincides with $f$. For that purpose set $V := {\pi'}^{-1} f^{-1} (\TT)$ and
consider the restriction $f \circ \pi'|_{V} \colon V \to \TT$. Since $f(X')
\cap \TT \neq \emptyset$ it suffices to show that this restriction coincides
with the corresponding restriction of $q \circ F_{t}$. Thus it is enough to
show that the following diagram of algebras of rational functions is
commutative:
\begin{equation}\label{kommtorus}
\xymatrix{%
{{\mathcal{R}} (\t{X'})} & {{\mathcal{R}} (\rq{X})} \ar[l]_{{F^{*}_{t}}} \\
{{\mathcal{R}} (X')} \ar[u]^{{\pi'}^{*}} & {{\mathcal{O}} (\TT) = \KK[M]}
\ar[l]^{f^{*}} \ar[u]_{q^{*}}.  }
\end{equation}
\begin{lemma}\label{texists}
There is a $t\in {(\KK^{*})}^{r}$ such that the diagram \eqref{kommtorus} is
commutative. Moreover, two such $t'$s differ only by an element of $H$.
\end{lemma}
\proof 
Let $m_{1},\dots, m_{n}$ denote a basis of $M$. Since $M \subset \widehat{M}$
each such $m_{i}$ can be considered as an element of $\widehat{M}$. More
precisely, using the basis $B_{1}, \dots, B_{r}$ introduced above, there is a
representation
\[ \div \chi^{m_{i}} = \sum_{j=1}^{r} a_{ij}B_{j}.\] 
An easy calculation shows that the map $q$ restricted to the torus of $\rq{X}$
is nothing but the map
\begin{equation}\label{quotientmap}
{(\KK^{*})}^{r} \to {(\KK^{*})}^{n}, \; t \mapsto \Bigl(\prod_{j=1}^{r}
t_{j}^{a_{1j}}, \dots, \prod_{j=1}^{r} t_{j}^{a_{nj}}\Bigr).
\end{equation}
As $q^{*}$ is the dual map of $q$ we explicitely can calculate $F_{t}^{*}
q^{*} (\chi^{m_{i}})$ and obtain
\begin{equation}\label{eq:normierung}
F_{t}^{*} q^{*} (\chi^{m_{i}}) = F_{t}^{*} \Bigl( \prod_{j=1}^{r}
{(T^{B_{j}})}^{a_{ij}}\Bigr) = \prod_{j=1}^{r} (t_{j}r_{j})^{a_{ij}} =
\prod_{j=1}^{r} t_{j}^{a_{ij}} F_{\boldsymbol{1}}^{*}q^{*}(\chi ^{m_{i}}),
\end{equation}
where $\boldsymbol{1} := (1, \dots, 1) \in ({\KK^{*}})^{r}$. On the other
hand, for $1\le i \le n$ by the definition of $F_{t}^{*}$ we obtain
\begin{equation}\label{eq:ident}
\div ({\pi'}^{*}f^{*} (\chi^{m_{i}})) = {\pi'}^{*} \phi (\div (\chi^{m_{i}}))
= F_{t}^{*} q^{*} (\div(\chi^{m_{i}})) = \div (F_{t}^{*} q^{*}
(\chi^{m_{i}})).
\end{equation}
For $t = \boldsymbol{1}$ this yields that the two functions ${\pi'}^{*} f^{*}
(\chi^{m_{i}})$ and $F_{{\boldsymbol{1}}}^{*} q^{*}(\chi^{m_{i}})$ differ only
by a constant, say $\lambda_{i}$. Thus, equation \eqref{eq:normierung} tells
us that, for $t \in {(\KK^{*})}^{r}$, \eqref{kommtorus} is commutative, if and
only if $\lambda_{i} = \prod_{j=1}^{r} t_{j}^{a_{ij}}$ for all $i$. Comparing
this with \eqref{quotientmap} this means that there is a $t \in
{(\KK^{*})}^{r}$ such that $q(t) = (\lambda_{1},\dots, \lambda_{n})$. But $q$
is the quotient map, which is geometric over the torus of $X$; therefore
$q|_{({\KK^{*})}^{r}}$ is surjective with kernel $H$. Hence the claim of the
Lemma follows.  \endproof 

We now finally give the proof of Theorem \ref{mainthm}.

\proof[Proof of Theorem \ref{mainthm}] 
We have already seen that liftings induce geometric pullbacks. Thus it remains
to show how to assign a lifting to a geometric pullback. Let $\phi \colon
\rq{M} \to \WD(X')$ be a geometric pullback with respect to $f$ such that
$\phi (\widehat{M}) \subset \widehat{M}' + \PDiv(X')$ holds. Define $F_{t}$ as
in \eqref{eq:ratpb} and choose according to Lemma \ref{texists} a $t \in
{(\KK^{*})}^{r}$ such that for $F := F_{t}$ the diagram \eqref{kommtorus} is
commutative. Then the equation $q \circ F = f \circ \pi'$ holds on $V =
{\pi'}^{-1} f^{-1} (\TT)$, and for continuity reasons also on
$\overline{V}$. But, since $f(X') \cap \TT \neq \emptyset$, the set $V$ is a
non-empty open set in $\t{X'}$; thus $\overline{V} = \t{X'}$, and $F$ is
indeed a lifting of $f$. Moreover, Lemma \ref{morprops}(ii) yields that $F$ is
equivariant with respect to $\psi$. Lemma \ref{texists} tells us that this
equivariant lifting is unique up to an $h \in H$. Thus we have defined
precisely one equivalence class of equivariant liftings $\t{X'} \to \rq{X}$.

From Corollary \ref{universal} we know that $\pi' \colon \t{X'} \to X$ factors
through $q' \colon \rq{X'} \to X'$ and $\pi'_{1} \colon \t{X'} \to
\rq{X'}$. Here $\pi_{1}'$ is a quotient presentation with respect to the
canonical operation of the group $G_{1} := \Hom (\Div_{W}^{\TT'} (X') /
\widehat{M}', \KK^{*})$ given by
\[ G_{1} \times {\mathcal{O}}(\t{X'}) \to {\mathcal{O}} (\t{X'}),\;
(g,{T}^{E}) \mapsto g(\overline{E})^{-1} {T}^{E}, \] where $\overline{E}$
denotes the class of $E$ in $\Div_{W}^{\TT'}(X') / \widehat{M}'$. Lemma
\ref{comorprops}(ii) gives that the functions $r_{j}$ are of the form
${T}^{D'_{j}}(s_{j} \circ \pi')$ for certain rational functions $s_{j}$ on
$X'$ and some divisors $D'_{j}$ in $\rq{M}'$. In particular,
$\overline{D'_{j}} = 0$ in $\Div_{W}^{\TT'}(X') / \widehat{M}'$. Thus, all
functions $r_{j}$ are invariant under the action of $G_{1}$ and therefore $F$
is $G_{1}$-invariant. Hence, $F$ induces a morphism $\widehat{F}$ with
$\widehat{F} \circ \pi_{1}' = F$. From the resulting commutative diagram
\[ \xymatrix{%
{\t{X'}} \ar[drr]^{F} \ar[dr]_{{\pi_{1}'}} \ar[dd]_{{\pi'}} \\ & {\rq{X'}}
\ar[r]_{\rq{F}} \ar[dl]_{q'} & {\rq{X}} \ar[d]^{q} \\ X' \ar[rr]^{f} & & X }
\] it follows that $\widehat{F}$ is indeed the desired lifting. The
equivariance of $\widehat{F}$ follows from the fact that $G_{1}$ is contained
in the kernel of the group homomorphism $\psi$ defined in Lemma
\ref{morprops}(ii).
\endproof
\section{Application to the case of Isomorphisms}
Recall that a morphism $X' \to X$ between toric varieties is called toric if
its restriction to the embedded tori induces a homomorphism of algebraic
groups. In this section we prove
\begin{thm}\label{algimpltor}
If two (not necessarily non-degenerate) toric varieties are isomorphic as
algebraic varieties, then they are also isomorphic as toric varieties.
\end{thm}

Assume firstly that $X$ and $X'$ are non-degenerate toric varieties. Let $\pi'
\colon \t{X'} \to X'$ and $\pi \colon \t{X} \to X$ denote the Cox's
constructions of $X'$ and $X$. The acting groups are denoted by $G'$ and $G$,
respectively. Then there is a nice characterisations of isomorphisms in terms
of liftings.

\begin{lemma}\label{isochar}
Let $f \colon X' \to X$ be a morphism, and assume that $F \colon \t{X'} \to
\t{X}$ is a lifting of $f$ which is equivariant with respect to $\psi \colon
G' \to G$. Then $f$ is an isomorphism if and only if $F$ and $\psi$ are
isomorphisms.
\end{lemma}
\proof {\bf ``only if:''} If $f$ is an isomorphism, then there exists a unique
geometric pullback $\phi$ with respect to $f$ (cf. Remark \ref{codim}). Since
$\phi$ is the strict transform of Weil divisors it is easy to see that $\phi$
is an isomorphism. Obviously, the induced map of divisor classes $[\phi]
\colon \ClDiv(X') \to \ClDiv(X)$ is also an isomorphism. Hence $\psi$ which is
the dual map to $[\phi]$ has to be an isomorphism, too. To show that $F$ is
also an isomorphism, let $F'$ denote a lifting of $f^{-1}$, which exists
according to Corollary \ref{isolift}. Then we have
\[\pi F F' = f \pi' F' = f f^{-1} \pi = \pi ,\] whence $FF'$ is a lifting of
the identity on $X$. But the identity admits upto equivalence a unique
lifting, thus $FF'$ differs only by a constant from the identity on
$\t{X}$. Therefore $F$ has a right inverse. By interchanging the roles of $F$
and $F'$ we obtain that $F$ also has a left inverse and we are done.

{\bf ``if:''} From the assumption that $F$ and $\psi$ are isomorphisms, it
follows immediately that $F^{-1}$ is equivariant with respect to
$\psi^{-1}$. Thus $F^{-1}$ defines a morphism $f' \colon X \to X'$. Similarly
as above
\[ f'f \pi' = f' \pi F = \pi' F^{-1}F = \pi'.\]
As $\pi'$ is surjective this yields $f'f = \id_{X'}$. In the same way one
proves $ff'= \id_{X}$ which establishes the assertion.  \endproof
Using the results obtained so far we are able to give a short proof of Theorem
\ref{algimpltor}.

\proof[Proof of Theorem \ref{algimpltor}] According to Lemma \ref{nondeg}
below we can assume that $X$ and $X'$ are non-degenerate toric varieties. Let
$f \colon X' \to X$ be an isomorphism. Due to Corollary \ref{isolift} there is
a lifting $F\colon \t{X'} \to \t{X}$, which is equivariant with respect to
some homomorphism $\psi$ of algebraic groups. From Lemma \ref{isochar} we know
that $F$ and $\psi$ are isomorphisms.

As the codimensions of the exceptional sets $Z$ and $Z'$ in the affine spaces
$\KK^{k}$ and $\KK^{k'}$, respectively, are at least two, $F$ extends to an
isomorphism $F \colon \KK^{k'} \to \KK^{k}$. In particular, $k = k'$. Since
$F$ is an isomorphism, the jacobian $J_{F} (0) = \det\Bigl(\frac{\partial
F}{\partial z} (0) \Bigr)$ cannot equal zero. Hence, there is a permutation $p
\in \text{Perm}\{1,\dots, k\}$ such that in each polynomial $F_{j} :=
\pr_{j}\circ F$ the monomial $z_{p(j)}$ occurs.

Obviously, $F$ induces an isomorphism of the exceptional sets $Z'$ and $Z$. In
particular, it maps each irreducible component of $Z'$ onto an irreducible
component of $Z$. Recall that an irreducible component of $Z'$ is of the form
\[ Z'(i_{1},\dots, i_{r}) = V(\KK^{k}; z'_{i_{1}}, \dots,
z'_{i_{r}}) = \bigcap_{\nu = 1}^{r} \{z'_{i_{\nu}} = 0\} \] where $r$ is
minimal such that the corresponding divisors $D_{i_{l}}$, $l =1, \dots, r$,
have an empty intersection.

We claim: If $F (Z'(i_{1},\dots, i_{r})) = Z(j_{1},\dots, j_{s})$, then $r =
s$ and $p (\{j_{1},\dots,j_{s}\}) = \{ i_{1},\dots, i_{r}\}$.

Proof: Since $\dim Z'(i_{1},\dots, i_{r}) = k-r$ and $\dim Z(j_{1},\dots,
j_{s}) = k-s $ it is obvious that $r=s$. \\
Assume now that there exists a $j_{*} \in \{j_{1},\dots,j_{s}\}$ with
$p(j_{*}) \notin \{ i_{1},\dots, i_{s}\}$. Then for the set $W := \bigcap_{i
\neq p(j_{*})} \{ z'_{i} = 0 \}$ the inclusion $W \subset Z'(i_{1}, \dots,
i_{s})$ holds. Applying $F$ yields
\[ F(W) \subset F(Z'(i_{1}, \dots, i_{s})) =Z(j_{1}, \dots, j_{s}) =
\bigcap_{\mu = 1}^{s} \{ z_{j_{\mu}}=0 \}, \] in particular this means
$F_{j_{*}} (W) = 0$. But as observed above the polynom $F_{j_{*}}$ is of the
form $F_{j_{*}}(z) = a_{0} + a_{1}z_{p(j_{*})} + R(z)$, where $z_{p(j_{*})}$
does not occur as a monomial in $R(z)$. Thus $F_{j_{*}}|_{W}$ is a
non-constant polynomial in $z_{p(j_{*})}$ contrary to $F_{j_{*}}(W) = 0$.

Set $F_{p} (z_{1},\dots, z_{k}) := (z_{p(1)}, \dots, z_{p(k)})$. The last
claim yields that $F_{p}$ induces a morphism $F_{p} \colon \t{X'} \to
\t{X}$. Moreover, since $F$ is equivariant with respect to $\psi$, all $F_{j}$
are homogeneous polynomials in the grading defined by the action of $G'$. As
$z_{p(j)}$ occurs in $F_{j}$, the degree of $z_{p(j)}$ with respect to this
grading coincides with the degree of $F_{j}$. But this implies that $F_{p}$ is
also equivariant with respect to $\psi$. Thus $F_{p}$ defines a morphism
$f_{p} \colon X' \to X$. Using Lemma \ref{isochar} the fact that both $F_{p}$
and $\psi$ are isomorphisms yields that $f_{p}$ also is an
isomorphism. Moreover, it is evident, that $F_{p}$ is a toric morphism. Since
$\pi'$ and $\pi'$ are even toric quotients, this eventually shows that $f_{p}$
is a toric isomorphism from $X'$ onto $X$.  \endproof
We are left to show the corresponding statement for degenerate toric
varieties. But for that purpose it suffices to prove the following
cancellation statement:
\begin{lemma}\label{nondeg}
Let $X'$ and $X$ be non-degenerate toric varieties, and assume that there is
an isomorphism $f \colon X' \times {(\KK^{*})}^{r'} \to X \times
{(\KK^{*})}^{r}$. Then we have $r = r'$ and $X \cong X'$ as algebraic
varieties.
\end{lemma}
\proof
For $u \in X'$, consider the morphism
\[ h_{u} \colon {(\KK^{*})}^{r'} \cong \{u\} \times
{(\KK^{*})}^{r'} \hookrightarrow X' \times {(\KK^{*})}^{r'}
\overset{f}{\longrightarrow} X \times {(\KK^{*})}^{r}
\overset{\pr_{2}}{\longrightarrow} {(\KK^{*})}^{r}.\] On the other hand, a
fixed $t \in {(\KK^{*})}^{r}$ determines a morphism
\[ g_{t} \colon X' \cong X' \times \{t\} \hookrightarrow X' \times
{(\KK^{*})}^{r} \overset{f}{\longrightarrow} X \times {(\KK^{*})}^{r}
\overset{\pr_{2}}{\longrightarrow} {(\KK^{*})}^{r}. \] As $X'$ is
non-degenerate, it is easy to see that the morphism $g_{t}$ has to be
constant. But this yields
\[ h_{u} (t) = g_{t} (u) = g_{t} (v) = h_{v} (t), \] that is, $h_{u} = h_{v}$
for all $u,v \in X'$. In particular, each $h_{u}$ is surjective (since $f$ is
surjective). This implies $r' \ge r$. Analogously, by considering $f^{-1}$,
one obtains $r \ge r'$, whence $r = r'$ and $\dim X = \dim X'$ holds. By
multiplying $h_{u}$ with a constant we can assume that $h_{u}$ is already a
group morphism. As a surjective homomorphism of algebraic tori of equal
dimension, $h_{u}$ has finite kernel, that is, the unit element
$\boldsymbol{1} \in ({\KK^{*}})^{r}$ has at most finitely many preimages, say
$t_{1},\dots, t_{l}$. Since $X \times \{\boldsymbol{1}\}$ is irreducibel, we
conclude $l = 1$ and so
\[ X \cong X \times \{\boldsymbol{1}\} = f(X' \times \{t_{1}\}) \cong X'
\times \{t_{1}\} \cong X',\] which eventually gives the assertion.
\endproof


\thebibliography{99999}

\bibitem[AHS]{Ajs} A. A' Campo, J. Hausen, S. Schr\"oer, {\em Quotient
Presentations and Homogeneous Coordinates for Toric Varieties},
math.AG/0005086, to appear in Math. Nachr.

\bibitem[Be]{Be} F. Berchtold, {\em Morphismen zwischen torischen
Variet\"aten}, Dissertation an der Universit\"at Konstanz,
Hartung-Gorre-Verlag, Konstanz (2001)

\bibitem[Bo]{Bo} A. Borel, {\em Linear algebraic groups}, second enlarged
edition, Springer GTM (1991)

\bibitem[Br-Ve]{BV} M. Brion, M. Vergne, {\em An Equivariant Riemann-Roch
Theorem for Complete Simplicial Toric Varieties}, J. Reine Angew. Math. {\bf
482}, 67-92 (1997)

\bibitem[Co-1]{Cox1} D. A. Cox, {\em The homogeneous coordinate ring of a
toric variety}, Jour. of. Alg. Geom. {\bf 4}, 17-50 (1995)

\bibitem[Co-2]{Cox2} D. A. Cox, {\em The functor of a smooth toric variety},
T\^ohoku Math. J. {\bf 47}, 251-262 (1995)

\bibitem[Da]{Da} V. I. Danilov, {\em Geometry of Toric Varieties}, Russian
Math. Surveys {\bf 33}, 97-154 (1978)

\bibitem[Dem]{Dem} M. Demazure, {\em Sous-groupes alg\'ebriques de Rang
Maximum du Groupe de Cremona}, Ann. sci. \'Ec. Norm. Sup., {\bf 4} 3, 507-588
(1970)

\bibitem[De]{De} A. S. Demushkin, {\em Combinatoric Invariance of Toric
Singularities}, Mosc. Univ. Math. Bull. {\bf 37, No. 2}, 104-111 (1982)

\bibitem[Fi]{Fi} J. Fine, {\em On Varieties isomorphic in codimension one to
Torus Embeddings}, Duke Math. J. {\bf 58} No.1, 79-88 (1989)

\bibitem[Ful]{Ful} W. Fulton, {\em Introduction to Toric Varieties}, Princeton
University Press (1993)

\bibitem[Gu]{Gu} J. Gubeladze, {\em The isomorphism problem for commutative
monoid rings}, J. Pure Appl. Algebra {\bf 129, No. 1}, 35-65 (1998)

\bibitem[Ka]{Ka} T. Kajiwara, {\em The functor of a toric variety with enough
invariant effective Cartier divisors}, T\^ohoku Math. J. {\bf 50}, 139-157
(1998)

\bibitem[Kr]{Kr} H. Kraft, {\em Geometrische Methoden in der
Invariantentheorie}, Vieweg Braunschweig (1984)

\bibitem[Od]{Oda} T. Oda, {\em Convex Bodies and Algebraic Geometry},
Springer-Verlag (1988)

\bibitem[O-M]{OM} T. Oda, {\em Lectures on Torus Embeddings and Applications},
Tata Inst. of Fund. Research {\bf 58}, Springer-Verlag, Berlin (1978)
\end{document}